\documentclass[10pt]{amsart}

\usepackage{amssymb}
\usepackage{amsfonts} 
\usepackage{latexsym}   
\usepackage{amssymb}    
\usepackage{amsmath}    
\usepackage{amsbsy}
\usepackage{amsgen}
\usepackage{amsfonts}
\usepackage{array}
\usepackage{wrapfig}

\usepackage{epsfig}
\usepackage{psfrag}
\usepackage[all]{xy}
\usepackage{hyperref}
\usepackage{subfig}

\usepackage{xcolor}
\usepackage{amssymb}
\usepackage{mathabx}
\usepackage{amsfonts}
\usepackage{hyperref}
\usepackage{amsmath}
\usepackage{graphicx}
\usepackage[all]{xy}%
\providecommand{\U}[1]{\protect\rule{.1in}{.1in}}

\newcommand{\R}{\mathbb{R}}
\newcommand{\dist}{\mathrm{dist}}

\theoremstyle{plain}

\newtheorem{algorithm}{Algorithm}[section]

\newtheorem{example}[algorithm]{Example}

\newtheorem{question}[algorithm]{Question}
\newtheorem{problem}[algorithm]{Problem}

\newtheorem{theorem} [algorithm] {Theorem}

\newtheorem{theoremlet'}[thm]{Theorem$'$}

\newtheorem*{FT}{Frankel's Theorem}

\usepackage{amssymb}
\usepackage{amsfonts} 
\usepackage{latexsym}   
\usepackage{amssymb}    
\usepackage{amsmath}    
\usepackage{amsbsy}
\usepackage{amsgen}
\usepackage{amsfonts}
\usepackage{array}
\usepackage{scalerel}

\usepackage{epsfig}
\usepackage{psfrag}
\usepackage[all]{xy}
\usepackage{amssymb}
\usepackage{mathabx}
\newtheorem{remark}[algorithm]{Remark}

\newtheorem*{msrconj}{Maximal Symmetry Rank Conjecture}

\newtheorem*{1/2MSR}{Half Maximal Symmetry Rank Theorem}

\newtheorem*{observe*}{Observation}

\def\qqq{\mathbb{Q}}
\def\rrr{\mathbb{R}}
\def\ccc{\mathbb{C}}
\def\zzz{\mathbb{Z}}

\def\RP{\rrr\mathrm{P}}
\def\CP{\ccc\mathrm{P}}
\def\HP{\hh\mathrm{P}}

\def\hh{\mathbb{H}}

\DeclareMathOperator{\diam}{diam}

\DeclareMathOperator{\Isom}{Isom}

\def\bdm{\begin{displaymath}}
\def\edm{\end{displaymath}}
\def\beq{\begin{equation}}
\def\eeq{\end{equation}}
\def\bes{\begin{equation*}}
\def\ees{\end{equation*}}
\def\epcm{\end{picture}\end{center}\end{minipage}}
\def\bpcm{\begin{minipage}{80pt}\begin{center}\begin{picture}}

\def\t2{T^2}

\def\f4{F_4}
\def\g2{G_2}

\def\p2{\frac{\pi}{2}}
\def\dist{\textrm{dist}}

\def\rk{\textrm{rk}}

\def\dim{\textrm{dim}}

\def\symrk{\textrm{symrk}}

 \numberwithin{equation}{section}
  \numberwithin{figure}{section}

\DeclareMathOperator{\curv}{curv}

\DeclareMathOperator{\Ric}{Ric}
\newtheorem*{SP}{Symmetry Program}
\newtheorem*{ST}{Soul Theorem}
\newtheorem*{STA}{Soul Theorem for Alexandrov Spaces}
\newtheorem*{GBNT}{Betti Number Theorem}
\newtheorem*{SplittingT}{Splitting Theorem}
\newtheorem*{LHT}{Lichnerowicz-Hitchin Theorem}
\newtheorem*{BMT}{Bonnet-Myers Theorem}
\newtheorem*{Synge}{Synge's Theorem}

\usepackage[margin=1.0in,footskip=0.25in]{geometry}

\begin{document}
\newcommand{\comment}[1]{\vspace{5 mm}\par \noindent
\marginpar{\textsc{Note}}
\framebox{\begin{minipage}[c]{0.95 \textwidth}
#1 \end{minipage}}\vspace{5 mm}\par}

\title[Symmetries of Spaces with Lower Curvature Bounds]{Symmetries of Spaces with Lower Curvature Bounds}

\author[Searle]{Catherine Searle}

\address[Searle]{Department of Mathematics, Statistics, and Physics, Wichita State University, Wichita, Kansas}
\email{catherine.searle@wichita.edu}


\date{\today}

%
\maketitle


\vspace{-0.5 cm}
\section{Introduction}

Global Riemannian Geometry generalizes the classical  Euclidean, Spherical and Hyperbolic geometries. One of the major challenges in this area is to understand how 
local invariants such as curvature, that is, how much a space ``bends",  relate to global topological invariants such as fundamental group, itself a measure of how ``connected" a manifold is. While locally Riemannian manifolds  have the topology of Euclidean space, the geometry typically deviates from that of $\rrr^n$: curvature is the cause of this deviation. 

In this article our main focus is on closed manifolds with bounded lower {\em sectional} curvature. A simple way to understand a lower {\it sectional} curvature bound is via triangle comparisons. We say that a manifold has a lower curvature bound $\kappa$ if the angle sum of any geodesic triangle, that is, a triangle formed by shortest length curves, 
is larger than or equal to the angle sum of the corresponding triangle in  $M^2(\kappa)$, the $2$-dimensional model space with constant curvature $\kappa$. Thus, we say that a manifold has positive, zero or negative  curvature, that is, $\kappa >0$, $\kappa=0$, or $\kappa<0$, respectively, if the angle sum of a geodesic triangle,  is strictly greater than $>\pi$, equal to $\pi$, or strictly less than $\pi$, respectively. 
In Figure \ref{triangles} below, we see how a geodesic triangle looks in positive, zero, and negative curvature, that is for $\kappa>0, \kappa=0$, and $\kappa<0$, respectively.

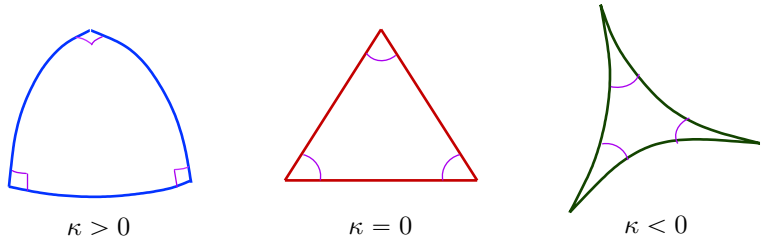
\begin{figure}[htbp]
\vspace{0.2cm}
\hspace{.7cm}
\begin{picture}(0,0)%
\includegraphics[scale=0.5]{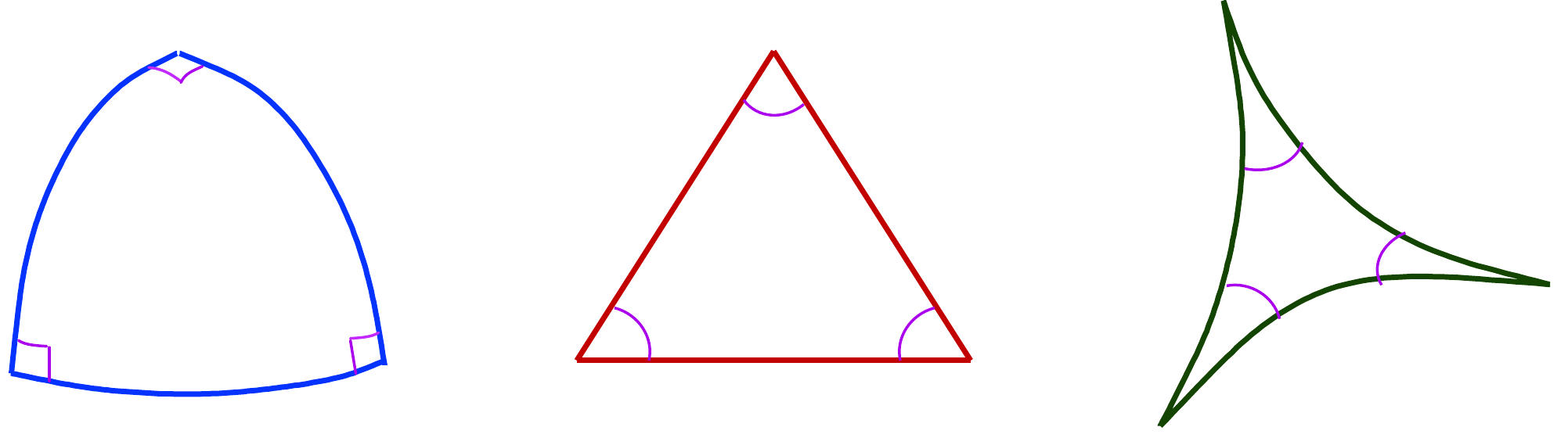}%
\end{picture}%
\setlength{\unitlength}{3947sp}%
\begingroup\makeatletter\ifx\SetFigFont\undefined%
\gdef\SetFigFont#1#2#3#4#5{%
  \reset@font\fontsize{#1}{#2pt}%
  \fontfamily{#3}\fontseries{#4}\fontshape{#5}%
  \selectfont}%
\fi\endgroup%
\begin{picture}(5079,1559)(1902,-7227)
\put(2300,-7336){\makebox(0,0)[lb]{\smash{{\SetFigFont{10}{8}{\rmdefault}{\mddefault}{\updefault}{$\kappa>0$}%
}}}}
\put(4070,-7325){\makebox(0,0)[lb]{\smash{{\SetFigFont{10}{8}{\rmdefault}{\mddefault}{\updefault}{$\kappa=0$}%
}}}}
\put(5800,-7320){\makebox(0,0)[lb]{\smash{{\SetFigFont{10}{8}{\rmdefault}{\mddefault}{\updefault}{$\kappa<0$}%
}}}}
\end{picture}%
\caption{Geodesic Triangles}
\label{triangles}
\end{figure}

\vspace{0.5cm}
\section{Closed Manifolds of Positive and Non-negative Sectional Curvature}

The classification of closed manifolds of positive and non-negative sectional curvature is a long-standing and very difficult problem in Riemannian geometry.
Unless otherwise stated, all curvatures considered here are sectional. For positive curvature, to date, other than some special examples in dimensions less than or equal to $24$, all known simply-connected examples are spherical in nature, that is, they are spheres, or quotients of spheres: $\mathrm{S}^n$, $\ccc\mathrm{P}^k$, or $\hh\mathrm{P}^l$.
There are many more examples of Riemannian manifolds of non-negative curvature. We list a few examples:
\begin{itemize}
\item Homogeneous spaces, $G/H$
\item Products of manifolds of non-negative curvature
\item Biquotients, $G//H$
\item Bases of Riemannian submersions of non-negatively curved manifolds.
\end{itemize}

When approaching classification problems, one looks for {\em obstructions} and {\em constructions}.
Among the obstructions to positive and non-negative curvature we have the following five, now classical, 
results.

\subsection{Obstructions}
We begin by listing obstructions for non-negative sectional curvature.  Note that some of these come from results about other types of lower curvature bounds such as Ricci, which is an average of the sectional curvatures, and scalar, which is an average of the Ricci curvatures. In particular, if $\sec(M)\geq \kappa$ then both the Ricci and scalar curvatures are bounded below by $\kappa$.

The first result is due to Cheeger and Gromoll  \cite{CG1}
 and tells us that when studying manifolds of non-negative curvature, it suffices to limit our attention to those that are closed, that is compact without boundary.
\begin{ST}[Cheeger and Gromoll \cite{CG1}] Let $M$ be a complete manifold of non-negative sectional curvature. Then there exists a closed,  totally geodesic, embedded submanifold, $S$, the {\em soul} of $M$, and $M$ is diffeomorphic to the normal bundle over $S$.
\end{ST}
\noindent The next theorem limits the fundamental group of a complete, non-negatively curved manifold and is also due to Cheeger and Gromoll \cite{CG2}.
 It was originally stated for compact manifolds of non-negative Ricci curvature.
\begin{SplittingT}[Cheeger and Gromoll \cite{CG2}] Let $M$ be a closed manifold of non-negative sectional curvature. Then the universal cover of $M$, $\widetilde M$, splits isometrically as 
$$\widetilde M=N\times \rrr^k,$$
where $N$ is a closed, simply-connected, non-negatively curved Riemannian manifold, and $k$ is the abelian rank of $\pi_1(M)$.
Moreover, $\pi_1(M)$ contains an abelian subgroup of finite index.
\end{SplittingT}
\noindent The next theorem, due to Gromov \cite{G}, 
limits the total Betti number of a manifold of non-negative curvature. As one application, it tells us that we can only take a limited number of connected sums of complex projective spaces and maintain non-negative sectional curvature.
\begin{GBNT} Let $M^n$ be a compact $n$-dimensional manifold of non-negative sectional curvature. Then  there exists a constant $\mathcal{C}(n)$ such that for any field $\mathbb{F}$,
$$\Sigma b_i(M^n; \mathbb{F}) \leq \mathcal{C}(n).$$
\end{GBNT}
\noindent We say a manifold is flat if all of its sectional curvatures are identically zero.
It is known that compact, flat, spin manifolds all have vanishing $\hat{A}$-genus and $\alpha$-invariant.
Thus, the topological obstructions for positive scalar curvature due to Hitchin and Lichnerowicz give us an obstruction for manifolds of non-negative curvature. 
\begin{LHT}  A compact, spin, $n$-dimensional manifold, $M^n$, with $\hat{A}(M)\neq 0$ or $\alpha(M)\neq 0$ does not admit a metric of non-negative sectional curvature. 
\end{LHT}
\noindent  For example, via this theorem, there are $9$-dimensional exotic spheres that carry no metric of non-negative sectional curvature. 

We now pass to obstructions for strictly positive sectional curvature.
The first is due to Bonnet 
and Myers, which was originally stated for manifolds with uniformly bounded positive Ricci curvature, and tells us that when studying manifolds of positive curvature, we may restrict our attention to those that are compact and simply-connected. 

\begin{BMT} Let $M$ be a Riemannian manifold with $\sec(M)\geq \kappa >0$, where $\kappa>0$. Then $M$ is compact and $\pi_1(M)$ is finite.
\end{BMT}
\noindent Finally, a result due to Synge gives us information about the fundamental groups and orientability of closed manifolds of positive curvature.

\begin{Synge}\label{SyT} Let $M$ be a closed manifold of positive sectional curvature. Then 
the following hold:
\begin{enumerate}
\item If $M$ is even-dimensional, then $\pi_1(M)$ is either trivial or $\zzz_2$; and 
\item If $M$ is odd-dimensional, then $M$ is orientable.
\end{enumerate}
\end{Synge}
\noindent In particular, Synge's Theorem tells us that $\RP^2\times \RP^2$ does not admit a metric of positive sectional curvature. Likewise, for a closed, orientable, odd-dimensional manifold $M^{2n+1}$, the product manifold, $M^{2n+1}\times \RP^2$ does not admit positive sectional curvature.

It bears mentioning at this point that despite the difference in  magnitude of the number of examples of positive and non-negative sectional curvature, when we restrict our attention to the class of closed, simply-connected manifolds, there are no manifolds that admit a metric of non-negative sectional curvature that are known to not admit a metric of positive sectional curvature. 

\subsection{Constructions}
Turning our attention to constructions, we have the Gray--O'Neill curvature equations and Cheeger deformations.
The Gray-O'Neill curvature equations tell us that for a Riemannian submersion, $\pi:E\rightarrow B$, if $\sec(E)\geq\kappa$ for some real number $\kappa$, then $\sec(B)\geq\kappa$. That is, curvatures can only increase. For example, via the Gray--O'Neill equations, it follows that any homogeneous space, $G/H$, admits a submersion metric of non-negative curvature from the bi-invariant metric $g_{\textrm{bi}}$ on the compact Lie group $G$. The same can be seen to be true for bi-quotients, $G//H$, defined as quotients of the free action of $H\subset G\times G$ on $G$ given by $(h_1, h_2)*(g)=h_1gh_2^{-1}$. Observe that while $S^2\times S^2$ covers $\RP^2\times \RP^2$, if $S^2\times S^2$ were to admit positive sectional curvature, see the first Hopf Conjecture below, then the covering map is {\em not} a Riemannian submersion.

\begin{figure}[!htb]
 \centering
 \subfloat[$l$ close to infinity]
 {%
      \includegraphics[width=0.25\textwidth]{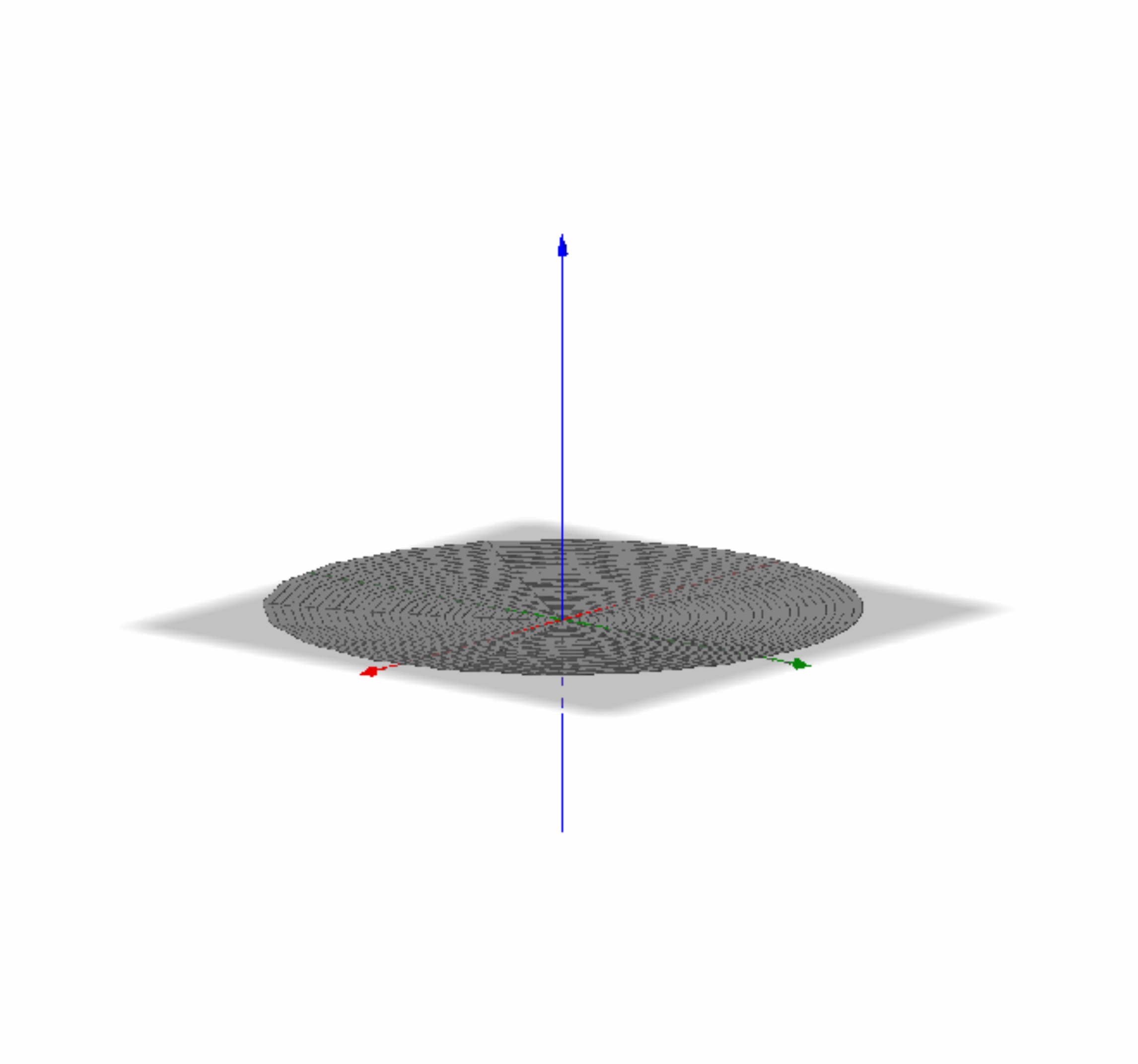}}
      \label{fig:image-a}
 \qquad
 \subfloat[$0<l<\infty$]
 {%
      \includegraphics[width=0.25\textwidth]{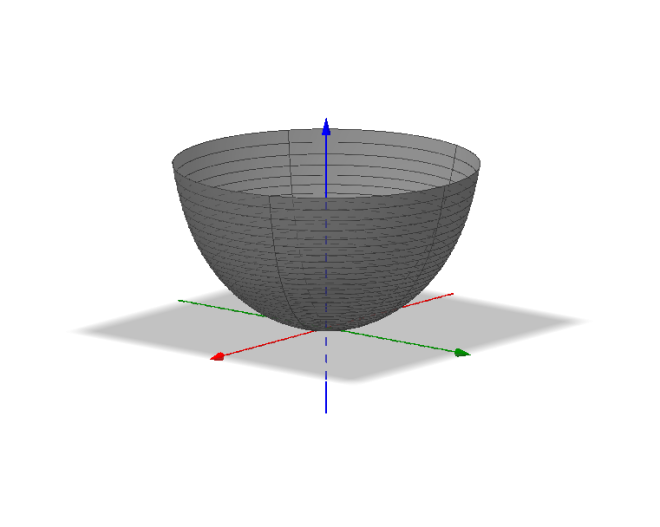}}
      \label{fig:image-b}
 \qquad
 \subfloat[$l$ close to $0$]
{%
      \includegraphics[width=0.25\textwidth]{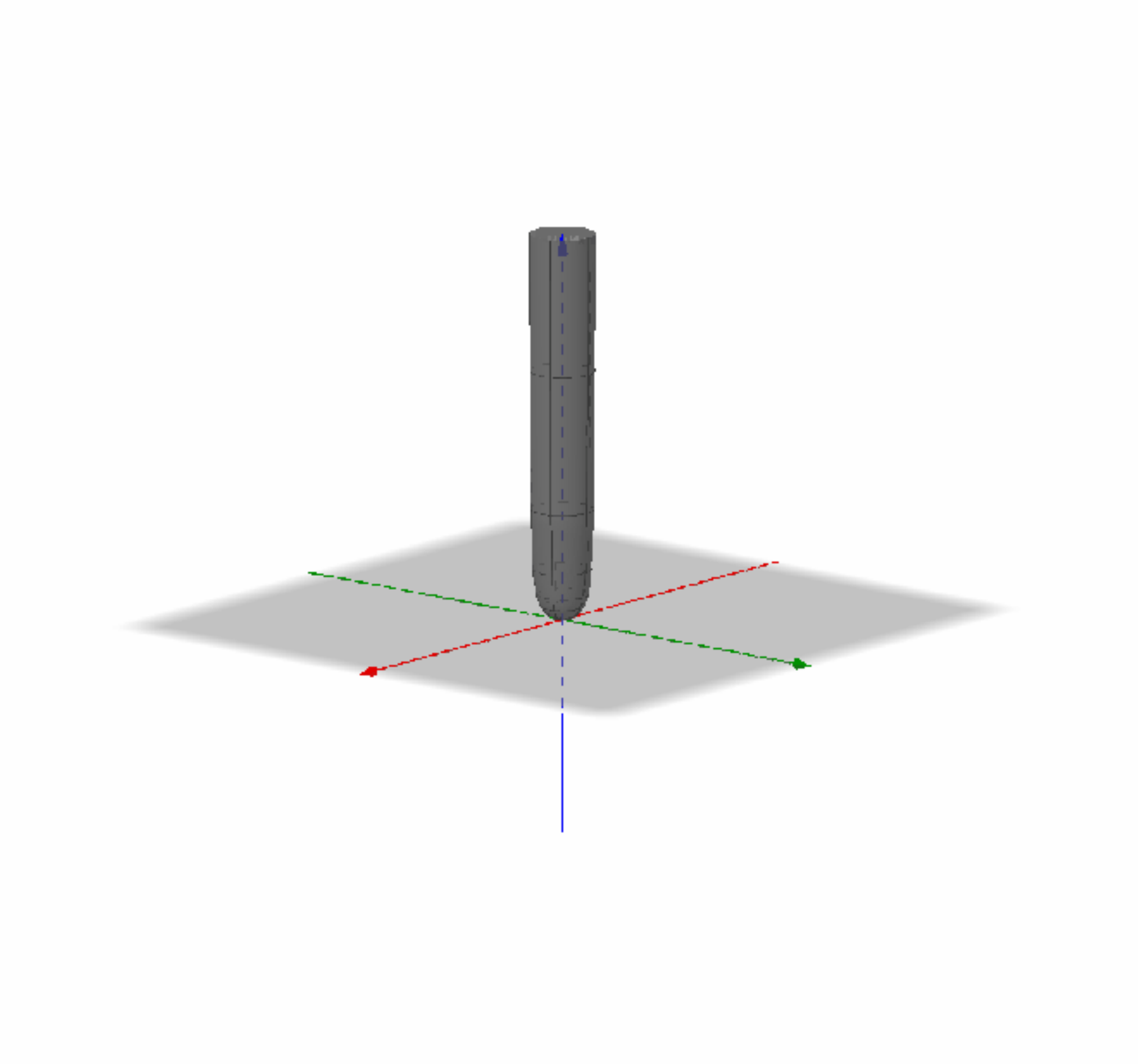}}
      \label{fig:image-c}
      \caption{A Cheeger deformation of the plane}
       \label{fig:Cheeger}
\end{figure}

Let $G$ be a closed subgroup of the isometry group of $M$, a Riemannian manifold,  endowed with any bi-invariant metric $g_G$. Cheeger deformations on a Riemannian manifold, $M$, leverage the power of Riemannian submersions 
by submersing from a $G$-extension of the manifold, $\pi:M\times G\rightarrow M$, where the base space of the Riemannian submersion is  
$M=(M\times G)/G$,
with the free $G$-action given by $g'(x, g)=(g'x, g'g)$. 
We obtain a one-parameter family of metrics on the total space, $\{(M\times G, g_M + l^2g_G)\}$, giving us a one-parameter family of metrics on the quotient space $\{((M\times G)/G, g_l)\}$. This family of metrics is called a Cheeger deformation of the original manifold $(M, g_M)$. As $l$ approaches infinity, the deformed metrics $g_l$ converge to the original metric $g_M$, and as $l$ approaches $0$, the sequence of manifolds $\{((M\times G)/G, g_l)\}$ converge to the quotient space $M/G$. Note that positive curvature is preserved under such deformations, while non-negative curvature may be improved to positive curvature.
\begin{example} Consider the linear circle action on $\ccc$ given by complex multiplication, where $\ccc$ has the flat metric. 
As $l$ decreases, the manifold admits strictly positive sectional curvature. Notably, in the limit, the manifold collapses to a half line. See Figure \ref{fig:Cheeger} below for three images describing this deformation\footnote{Courtesy of Lawrence Mouill\'e, see \url{https://lawrencemouille.wordpress.com/author/lawrencemouille/} for a .gif of this action}.

\end{example}

Cheeger deformations have been useful in deforming $G$-invariant metrics and were originally used  by Cheeger  \cite{Ch}
 to show that $\CP^2\#\CP^2$ admits a metric of non-negative sectional curvature.  In particular, while not named explicitly, they feature in the theorem of Lawson and Yau  \cite{LY} 
 showing that the existence of a smooth non-abelian compact Lie group action guarantees positive scalar curvature on a compact manifold. They are also an important tool for many lifting theorems, where one tries to lift a synthetic curvature lower bound on the quotient space $M/G$ to the manifold $M$, such as in work of Searle and Wilhelm \cite{SW}.  
Additionally, they have recently been used by Cavenaghi, Grama, and Speran\c{c}a \cite{CGSp} who claim to show that the base of a positively curved principal $SU(2)$- or $SO(3)$-bundle must have dimension greater than or equal to $4$. If the result is true,  it provides an answer to a special case of the Petersen-Wilhelm conjecture that states that for a fibration of a positively curved manifold, the dimension of the fiber must be strictly less than the dimension of the base.

\subsection{Important Open Questions}
Finally, three important  open questions for manifolds of positive and non-negative curvature are contained in the following conjectures.
\begin{itemize}
\item (H. Hopf) {\em $S^2\times S^2$ does not admit a metric of positive sectional curvature.}
\item (H. Hopf) {\em Let $M^{2m}$ be an even-dimensional, closed Riemannian manifold of positive, respectively, non-negative sectional curvature. Then $M^{2m}$ has positive, respectively, non-negative Euler characteristic.}
\item (Bott){\em Let $M^n$ be a closed, simply-connected Riemannian manifold of non-negative sectional curvature. Then $M$ is rationally elliptic.}
\end{itemize}
Recall that a closed, simply-connected topological space is called {\em rationally elliptic} if $\pi_*(X)\otimes \,\qqq$ and $H_*(X; \qqq)$ are finite-dimensional $\qqq$-vector spaces. If we drop the simply-connected hypothesis, we refer to the space as {\em rationally $\Omega$-elliptic}.  Rationally elliptic manifolds 
have been classified in dimensions less than or equal to $5$ and are diffeomorphic to the known examples of manifolds of non-negative curvature.
The classification in dimensions $4$ and $5$ is curvature-free and due to Paternain and Petean   \cite{PaPe}, 
while in dimensions less than or equal to $3$ we only have spheres by the Gauss-Bonnet Theorem in dimension $2$ and the work of Hamilton  \cite{Ha}
 in dimension $3$. 
Rationally elliptic manifolds have non-negative Euler characteristic and if the Euler characteristic is positive, then all odd degree Betti numbers vanish. Resolving the Bott Conjecture would then resolve the second Hopf Conjecture for non-negatively curved manifolds.

Note that in dimensions greater than or equal to $4$, the classification of closed, simply-connected, positively and non-negatively curved manifolds is still open.
In an attempt to address this issue in dimension $4$, we have the following theorem due to Hsiang and Kleiner \cite{HsKl}. 
\begin{theorem}[Hsiang and Kleiner] \cite{HsKl}
\label{HK} Let $M^4$ be a closed, simply-connected, $4$-manifold of positive sectional curvature admitting an isometric and effective $T^1$-action. Then $M$ is homeomorphic to $S^4$ or $\CP^2$.
\end{theorem}

The proof of this theorem reduces to proving that the Euler characteristic of the manifold is bounded between $2$ and $3$ and applying deep classification work of Freedman  \cite{Fr}. 
In particular, the theorem tells us that if $S^2\times S^2$ were to admit a metric of positive curvature, then it must have a finite group of isometries.

The theorem has since been improved to diffeomorphism by work of Grove and Searle  \cite{GrS1} 
and Grove and Wilking  \cite{GrWi}. 
It was extended to non-negative curvature by independent work of Kleiner \cite{Kl} 
and Searle and Yang \cite{SY}, 
Galaz-Garc\'ia  \cite{G-G}, 
Galaz-Garc\'ia and Kerin  \cite{G-GK}, 
and Grove and Wilking  \cite{GrWi}.  
There, one sees that three more manifolds occur: $S^2\times S^2$ and $\CP^2\#\pm \CP^2$. More recently, the theorem 
was extended to almost non-negative curvature by Harvey and Searle \cite{HaS}, who showed that only the same manifolds occur.  Grove and Halperin suggested extending the Bott Conjecture to include almost non-negatively curved manifolds.
The result in \cite{HaS} confirms this extended Bott Conjecture with the addition of $S^1$-symmetry in dimension $4$.

\subsection{The Symmetry Program}

In the early nineties, Karsten Grove,  inspired by Theorem \ref{HK}  and observing that the known examples at that time of positive and non-negative curvature  were all highly symmetric, proposed his {\em Symmetry Program}:

\begin{SP} Classify closed manifolds of positive and non-negative curvature with ``large" symmetries.
\end{SP}

By work of Myers and Steenrod \cite{MySt}, 
 the isometry group of a compact manifold is 
a compact Lie group, so when we talk about ``symmetries", we mean an isometric action by a compact Lie group.
An attractive aspect of this program is the flexibility of the term ``large".
For example, for a given $G$-action on a Riemannian manifold $M$,  ``large" can mean that the dimension of the quotient space, $M/G$, is small. Another perspective is to consider $G$-actions with large fixed point sets, and yet another is to consider $G$-actions with large rank.
The goal of this program is to successively lower the size of the group and in the process find new examples, new obstructions, or new constructions. To date the program has been quite successful in positive and non-negative curvature and has been extended to other types of lower curvature bounds, as well as other spaces that generalize Riemannian manifolds. We will survey some of the
results leading to partial classifications, as well as partial resolutions of the three conjectures mentioned earlier.

For further information on the subject of positive and non-negatively curved manifolds with symmetries,   there are surveys by  Grove \cite{Gr:S}, 
Ziller \cite{Z1}, Ziller \cite{Z2}, 
and Wilking \cite{Wi5}.

\section{Preliminaries}

Before we begin, we first establish some notation as well as some useful facts about transformation groups in general. Remark that manifolds are assumed to be closed unless otherwise specified.

\subsection{Transformation Groups}

Let $G$ be a compact Lie group acting on a smooth manifold $M$. We denote by $G_x=\{\, g\in G : gx=x\, \}$ the \emph{isotropy group} at $x\in M$ and by $G(x)=\{\, gx : g\in G\, \}\simeq G/G_x$ the \emph{orbit} of $x$.  Orbits are called {\em principal, exceptional,} or {\em singular} depending on the size of their isotropy group, as follows. An orbit is principal  if the isotropy group is the smallest possible among all isotropy groups. Orbits are called exceptional when their isotropy group is a non-trivial finite extension of the principal isotropy subgroup, and singular when their isotropy group is of strictly larger dimension than that of the principal orbits. The isotropy subgroups of an orbit are conjugate to one another, that is, given $y\in G(x)$, $y=gx$ for some $g\in G$, and $g^{-1}G_yg=G_x$. It makes sense then to talk of the {\em isotropy type} of an orbit. A $G$-action on a manifold defines a natural stratification, with strata corresponding to the union of orbits of each isotropy type.

For isometric actions of compact Lie groups, the Slice Theorem  gives us an explicit description of a sufficiently small tubular neighborhood of any orbit. Namely, given an orbit $G(x)\subset M$, a sufficiently small $\epsilon$-tubular neighborhood,  $D_{\epsilon}(G(x))$ 
is diffeomorphic to 
$G\times_{G_x} D^{\perp}_x,$ where $D^{\perp}_x$ is the corresponding $\epsilon$-ball at the origin of the normal space $T^{\perp}_x$ to $G(x)$ at the point $x$, called the normal slice to the orbit. 

For more details on the theory of transformation groups see Bredon \cite{Br}. 

\subsection{Alexandrov Spaces}
An Alexandrov space, $(X, \dist)$, is a finite dimensional length space with curvatures bounded below via triangle comparisons with the corresponding model spaces. 
All closed Riemannian manifolds are Alexandrov spaces and limits of Gromov-Hausdorff sequences of closed Riemannian manifolds with the same lower curvature bound are, also. 

For closed manifolds with sectional curvature bounded below by $\kappa$ and admitting an isometric $G$-action, the quotient space $M/G$ is an Alexandrov space with curvature bounded below by $\kappa$ (see Perelman and Petrunin \cite{PerPet})
 with locally totally geodesic orbit strata.  There is also a Soul Theorem for Alexandrov spaces  of  non-negative curvature due to Perelman \cite {Per}.
\begin{STA} Let $X^n$ be a complete, $n$-dimensional Alexandrov space with boundary $N$ with $\curv\geq 0$.  Then there exists a convex subset $S\subset X$, the {\em soul} of $X$, at maximal distance from $N$ and a deformation retraction of $X$ onto $S$.
\end{STA}

For more basic information about Alexandrov spaces, see Burago, Burago, and Ivanov \cite{BBI} 
and Alexander, Kapovitch, and Petrunin \cite{AKP}.

\subsection{Fixed Point Sets}
 We will denote by $M^G=\{\, x\in M : gx=x  \,\text{for all} \,g\in G \, \}$  the \emph{fixed point set} of the $G$-action. If $G=T^k$, then we will often simply write $M^T$ for its fixed point set.  
 
 The components of the fixed point set of an isometry are closed,  totally geodesic submanifolds of $M$ and orientable if $M$ is, by work of Kobayashi \cite{Ko}.  
 In the special case of a circle action, the components of $M^{S^1}$ are of even codimension. Note as well that for $H$ a proper subgroup of $G$, $M^G\subset M^H$ and  $N_G(H)$, the normalizer of $H$ in $G$, acts invariantly on $M^H$ with ineffective kernel $H$.  Thus, there is an induced action of $N_G(H)/H$ on $M^H$. If, moreover, $N_G(H)/H$ is connected, then the action is invariant on each connected component of $M^H$. So, for the special case where $G=T^k$ and $M$ is orientable, for every subtorus $T^l\subset T^k$, every $N\subset M^{T^l}$ is also an orientable, closed submanifold admitting an induced $T^{k-l}$-action with the same lower curvature bound and the same parity of dimension as $M$. This sets the stage for induction arguments, something quite unusual in Riemannian geometry.

Finally, we mention some results about the existence of fixed point sets in positive curvature.
A result of Weinstein   \cite{We} 
guarantees the existence of a fixed point for any orientation preserving isometry on an even-dimensional, orientable, closed manifold of positive curvature. This result was generalized to general torus actions in even dimensions by Berger \cite{B2}
 and in odd dimensions by Sugahara \cite{Su} 
 as follows.
\begin{theorem} Let $M^n$ be a closed, $n$-manifold of positive sectional curvature admitting an isometric, effective $T^k$-action. Then the following hold:
\begin{enumerate}
\item $($Berger$)$ \cite{B2} 
If $n=2m$, then $M^{T}\neq \emptyset$; and
\item $($Sugahara$)$ \cite{Su} 
If $n=2m+1$, then either there is a point $p\in M$ such that $T^k(p)\cong S^1$, or $M^{T}\neq \emptyset$.
\end{enumerate}
\end{theorem} 
\begin{remark}
The results of Berger and Sugahara do not hold for $\zzz_p^k$-actions. For example, letting $k=p=2$, there exists a $\zzz_2^2$ action on $S^2$ which has no fixed points, generated by the orientation-preserving involutions
$$\begin{pmatrix} -1 &&\cr
&-1&\cr
&&1
\end{pmatrix}\,\,{\textrm{ and }}\,\,
\begin{pmatrix} -1 &&\cr
&1&\cr
&& -1
\end{pmatrix}.$$ 
 Likewise, there are free $\zzz_2^2$ and $\zzz_3^2$ actions on closed, positively curved $7$-manifolds by work of Shankar   \cite{Sh} 
 and  Grove and Shankar \cite{GrSha}, 
 respectively  and free $\zzz_3^2$ actions on closed, positively curved $13$-dimensional manifolds by work of Bazaikin \cite{Ba}. 
\end{remark}

\begin{remark}\label{bottGGKR}
The results of Berger and Sugahara no longer hold when we pass to non-negative curvature: there exist free $T^k$ actions on the $k$-fold product of three-spheres. 
On the other hand, if the Bott Conjecture holds, then for an effective torus action of sufficiently large rank on a closed, simply-connected manifold of non-negative curvature, the torus action will have non-trivial isotropy. 

This is quantified in the 
work of Galaz-Garc\'ia, Kerin, and Radeschi \cite{G-GKR}, who show that
 if $M^n$, a  rationally elliptic $n$-dimensional smooth manifold, admits a smooth and effective $T^k$-action, then
 $k\leq \lfloor \frac{2n}{3}\rfloor$, and any subtorus acting freely on $M^n$ has rank bounded above by $\lfloor  \frac{n}{3}\rfloor$. 
 \end{remark}

\subsection{What is ``Large"?}

Here we give three examples of what  ``large" symmetries can mean. They are:
\begin{itemize}
\item {\em Small quotient space}, that is, $\dim(M/G)$ is small;
\item {\em Large fixed point set}, that is, $\dim(M^G)$ is large with respect to the dimension of the manifold $M$; and 
\item {\em Large rank}, that is, we consider group actions $G$ for which $\rk(G)$ is large with respect to the maximal possible rank of a group action on a manifold. 
\end{itemize}
In what follows, we will discuss these three different meanings of large and survey results for these definitions in both positive and non-negative curvature. 
We note that for the last definition, if one passes to discrete groups, for example, $\zzz_p^k$, we can define large  discrete $p$-symmetry rank, for $p$ a prime,  to be large with respect to the maximal possible 
 number $k$ such that  the isometry group of $M$ contains an elementary $p$-group of rank $k$.
 
\section{Small quotient space}
We survey results for $G$-actions on manifolds with small quotient spaces in both positive and non-negative curvature. Here, the general strategy is to leverage knowledge of the quotient space to identify the manifold.

\subsection{Positive Sectional Curvature}

This constraint has been utilized to obtain classifications of closed manifolds of positive sectional curvature of low cohomogeneity, where the cohomogeneity of a $G$-action on a manifold, $M$, is equal to the dimension of the quotient space, or equivalently, the codimension of the principal orbits, $G/H$ of the $G$-action. In particular, homogeneous spaces, those of cohomogeneity $0$, have been classified by Berger \cite{B1}, 
B\'erard-B\'ergery \cite{B-B}, 
Aloff and Wallach \cite{AW}, 
Wallach \cite{W}, 
Wilking \cite{Wi1}, 
and Wilking and Ziller \cite{WiZ}.

For cohomogeneity one, closed manifolds of positive sectional curvature have been classified in dimension $5$ by  
Searle \cite{Se}, 
in even dimensions by Verdiani \cite{V1, V2}, 
and in all odd dimensions but $7$ by Grove, Wilking, and Ziller \cite{GrWiZ}. 
Additionally, a list of possible candidates for dimension $7$ is given in Grove, Wilking, and Ziller. They are grouped into one isolated $7$-manifold, $R$, and two infinite families, $P_k$ and $Q_k$, with $P_1=S^7$ \cite{GrWiZ}. 
These $7$-dimensional candidates correspond to the total space of the Konishi bundle of the self dual Hitchin orbifold $\mathcal{O}_k$, see Hitchin \cite{H}. 
Of these candidates, a metric of positive curvature was found to exist on $P_2$, an exotic $T^1S^4$, homeomorphic, but not diffeomorphic, to the unit tangent bundle of the $4$-sphere,  independently by  Dearricott \cite{D}, 
and Grove, Verdiani, and Ziller \cite{GrVZ}. 
Verdiani and Ziller   \cite{VZ}
have also shown that $R$ does not admit a $G$-invariant cohomogeneity one metric of positive sectional curvature. More recently, Dearricott \cite{D2} claims to have shown that the remaining $P_k$ and all of the $Q_k$ do admit metrics of positive curvature, although none admit a $G$-invariant cohomogeneity one metric of positive curvature. 

When working with group actions, one may divide them into those that are {\em polar}, namely, those that admit a section, a closed, totally geodesic immersed submanifold that meets all orbits orthogonally, and those that are {\em non-polar}. Note that all transitive actions and all actions of cohomogeneity one are polar, but there exist cohomogeneity two actions on spheres that are non-polar.

\begin{example}  Viewing $\rrr^3$ as the set of self-adjoint $2 \times 2$ complex matrices, $$\rrr^3\cong V=\left\{X =\begin{pmatrix} x_1 & x_2 + ix_3 \\ x_2-ix_3 & -x_1 \end{pmatrix}; x_1, x_2, x_3 \in \rrr\right\},$$ we obtain an action of $U(2)$ action on $S^6\subset \rrr^3\oplus \ccc^2$, as follows.  
Let $A\in U(2)$ act by conjugation on $X\in V$  and by matrix multiplication on $(z_1, z_2)\in \ccc^2$. The quotient space $S^6/U(2)$ is homeomorphic to a $2$-disk, $D^2$, with a circle's worth of singular orbits with $T^1$-isotropy corresponding to $\partial D^2$ and an isolated vertex orbit with isotropy $T^2$ corresponding to a vertex on $\partial D^2$, see Bredon \cite{Br}. 
For more examples of non-polar cohomogeneity two actions on spheres, see work of Straume \cite{Str2, Str3}. 
\end{example}

In Fang, Grove, and Thorberggson \cite{FGrT}, they show that a   closed, simply-connected, Riemannian $n$-manifold $M^n$ of positive curvature with a polar $G$-action of cohomogeneity $\geq 2$, is equivariantly diffeomorphic to a compact rank one symmetric space (CROSS) with the corresponding linear $G$-action. Thus, in order to classify the remaining cohomogeneities up to $n-1$, it remains to consider those actions that are non-polar. 
While the quotient spaces of polar actions on manifolds of positive curvature necessarily have boundary by work of Wilking \cite{Wi3}, 
a result  later  generalized to singular Riemannian foliations by Corro and Moreno \cite{CoMo}, the quotient spaces of non-polar actions  may  or  may not have 
 boundary.

\begin{figure}[!htb]
 \centering
 \includegraphics[width=0.35\textwidth]{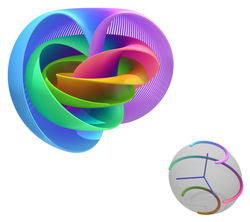}
       \caption{Stereographic Projection of The Hopf fibration of $S^3$ by $S^1$ and the Projection of its Fibers onto $S^2$}
        \label{Hopf}
\end{figure}

\begin{example} Some of the best known examples of positively curved manifolds are obtained as the base of a Riemannian submersion and correspond to the quotient space of a non-polar action, for example,  
$\CP^n$ and $\HP^k$ are both quotients of the Hopf action by $S^1$ and by $S^3$, respectively, on a sphere. Figure \ref{Hopf} depicts the Hopf fibration of $S^3$.\footnote{Image obtained from \url{https://en.wikipedia.org/w/index.php?title=Hopf_fibration&oldid=1101601722}}
\end{example}
A natural next step would be to tackle the following problem:
\begin{problem}
Classify  non-polar manifolds of positive curvature of low cohomogeneity, beginning with cohomogeneity two. 
\end{problem} 
\subsection{Non-negative sectional curvature}

As mentioned earlier, we can put a metric of non-negative sectional curvature on any homogeneous space $G/H$. In contrast, closed, simply-connected, cohomogeneity one manifolds of non-negative curvature have not yet been classified. They naturally admit a $G$-invariant disk bundle decomposition over each of the two singular orbits. While both disk bundles admit a $G$-invariant metric of non-negative sectional curvature,  in general, their union does not. Grove, Verdiani, Wilking, and Ziller   \cite{GrVWiZ} 
showed that some odd-dimensional exotic spheres, while they admit cohomogeneity one actions, do not admit such metrics. Later, C. He   \cite{He} 
showed that a larger class of manifolds that includes those in the work of Grove, Verdiani, Wilking, and Ziller   \cite{GrVWiZ} 
do not admit a cohomogeneity one metric of non-negative curvature.  
While closed, simply-connected,  manifolds admitting cohomogeneity one actions  may not admit invariant metrics of non-negative curvature, they do admit invariant metrics of almost non-negative curvature by work of Schwachh\"ofer and Tuschmann \cite{SchTu}.  
We say that a manifold, $M^n$,  is {\em almost non-negatively curved} if there exists a sequence of metrics, $\{g_n\}$ on $M$ and a fixed  $D>0$ so that 
	$	\diam(M,g_{n }) \leq D$ and 
		$\sec(M,g_n) \geq -\frac{1}{n^2}$.  
		Moving on to cohomogeneity two actions, we see already 
 in dimension $4$ that there are examples that do not admit non-negative curvature, such as $\CP^2\#\CP^2\#\CP^2$ with a $T^2$ isometric and effective action.

For a fibration, it is known that if any two elements of the fibration are rationally $\Omega$-elliptic, then so is the third.  Similarly, as observed by Grove and Halperin \cite{GrH}, 
for a manifold decomposing as a union of disk bundles 
$$M=D(A)\cup_E D(B),$$ if one of $A$, $B$, or $E$ is rationally $\Omega$-elliptic, (and hence are all, as the common boundary, $E$, is a sphere bundle over both $A$ and $B$) then $M$ is. Compact Lie groups are known to be rationally $\Omega$-elliptic. Moreover, a cohomogeneity one manifold decomposes as a union of disk bundles over its two singular orbits,  glued along the principal orbit. Thus  one sees that 
a closed, simply-connected manifold that is homogeneous or of cohomogeneity one is  rationally elliptic, regardless of curvature.

Grove and Halperin proposed that the Bott Conjecture will continue to hold for manifolds of almost non-negative curvature. 
 While classifying non-negatively curved manifolds of cohomogeneities $\geq 1$ seems out of reach at the moment, in light of the Bott Conjecture, asking a different question seems more tractable.
\begin{question}
Let $M$ be a  closed, simply-connected, almost non-negatively curved manifold. If $M$ admits an isometric $G$-action of  low cohomogeneity, is $M$ rationally elliptic?
\end{question}
This question has already been answered affirmatively for cohomogeneity two manifolds of almost non-negative curvature, by Grove, Wilking, and Yeager \cite{GrWiY}. Part of the proof of this result has been considerably simplified by recent work of  Khalili Samani and Radeschi  \cite{KSR} on {\em singular Riemannian foliations}, an area of study which can be viewed as a generalization of the concept of a group action.

\section{Large Fixed Point Set}

Here we survey results for $G$-actions on manifolds with large fixed point sets. Once again, one of the main strategies employed is to leverage an understanding of the quotient space to obtain general structure theorems that potentially lead to classification theorems.
\subsection{Positive Curvature}
An important first example of how large fixed point sets in positive curvature may limit the group action under consideration is given by Frankel's theorem, which tells us fixed point set components of sufficiently large dimension are unique in positive curvature.
\begin{FT} Let $M^n$ be an $n$-dimensional closed Riemannian manifold admitting a metric of positive sectional curvature. Suppose that $N^{k_1}_1$ and $N^{k_2}_2$ are two totally geodesic, embedded submanifolds of $M$. 
Then if $k_1+k_2\geq n$, $$N_1\cap N_2\neq \emptyset.$$
\end{FT}
We now observe that the dimension of the quotient space, $M/G$,  is constrained by the dimension of the fixed point set $M^G$  of $G$ in $M$. In fact, $\dim (M/G)\geq \dim(M^G) +1$ for any non-trivial, non-transitive action. In light of this, the {\it fixed-point cohomogeneity} of an action, denoted by $\textrm{cohomfix}(M;G)$, is defined by
\[
\textrm{cohomfix}(M; G) = \dim(M/G) - \dim(M^G) -1\geq 0.
\]
A manifold with fixed-point cohomogeneity $0$ is also called a {\it $G$-fixed point homogeneous manifold}.
 In Grove and Searle \cite{GrS2}, 
 they combine the critical point theory for distance functions introduced by Grove and Shiohama   \cite{GrSh} 
 (see also a survey on the subject by Grove \cite{Gr}) 
 and the Soul Theorem for Alexandrov spaces, to prove that for positively curved fixed point homogeneous manifolds there are at most three orbit types. These are given by the principal orbits, the fixed points contained in $F$, the fixed point set component of dimension equal to $\dim(M/G)-1$, and the unique orbit at maximal distance from the fixed point set component $F$. They use this to obtain a $G$-equivariant double disk bundle decomposition of the manifold, namely $M$ decomposes as 
$$M=D(F)\cup D(G(x)),$$
where $G(x)$ is the orbit in $M$  at maximal distance from $F$. They then use this decomposition to classify closed, simply-connected fixed point homogeneous manifolds of positive curvature, proving that for connected $G$, such a manifold is equivariantly diffeomorphic to a CROSS with a linear $G$-action.

This theorem was extended to the case of closed, simply-connected, fixed point cohomogeneity one manifolds of positive curvature by  Grove and Kim \cite{GrK}, 
who showed that such a manifold is also diffeomorphic to a CROSS. 
The fixed point homogeneous result has also been generalized to the case of involutions:  Fang and Grove \cite{FGr} showed that closed, positively curved, $\zzz_2$-fixed point homogeneous manifolds are diffeomorphic to spheres and real projective spaces.

\subsection{Non-negative curvature}

For the class of closed, simply-connected non-negatively curved manifolds, the fixed point homogeneous results were generalized to non-negative curvature in low dimensions by Galaz-Garc\'ia   \cite{G-G} 
and by 
Galaz-Garc\'ia  and Spindeler \cite{G-GSp}.
In his thesis, Spindeler    \cite{Spi} 
was then able to fully generalize the disk bundle decomposition result in the following theorem.
\begin{theorem}[Spindeler]\label{Spindeler}   \cite{Spi} 
Let $G$ act fixed point homogeneously  on a closed, non-negatively curved Riemannian manifold $M$. Let $F$ be a fixed point component of maximal dimension. Then there exists a smooth submanifold $N$ of $M$, without boundary, such that $M$ is diffeomorphic to the normal disk bundles $D(F)$ and $D(N)$ of $F$ and $N$ glued together along their common boundaries;
\bdm
M  = D(F) \cup_{E} D(N).
\edm
Further, $N$ is $G$-invariant and all points of $M \setminus \{F \cup N\}$ belong to principal $G$-orbits.
\end{theorem}

By contrast with the positively curved case, the submanifold $N$ at maximal distance from $F$ is not in general a single orbit.
   As noted by Spindeler, a classification of fixed point homogeneous manifolds of non-negative curvature in higher dimensions is currently out of reach since it is equivalent to a classification of non-negatively curved manifolds.  To see this, observe that for  $N$, a closed, simply-connected non-negatively curved manifold, the product manifold $M^{n}=N^{n-2}\times S^2$ with the product metric admits an isometric $S^1$-fixed point homogeneous action. Since we only have a classification of closed, simply-connected non-negatively curved manifolds through dimension $3$, obtaining a classification in dimension $5$ is the best we can hope for at the moment.

 \section{Large Symmetry and Discrete Symmetry Rank}

 In this section we discuss two notions of large rank. The first is large {\em symmetry rank}, where the symmetry rank of a $G$-action on $M$ is defined to be the rank of the isometry group of $M$, that is 
$$\symrk(M)=\rk(\Isom(M)).$$ 
The second is large {\em discrete $p$-symmetry rank}, for $p$ a prime, defined to be 
the largest number $k$ such that  the isometry group of $M$ contains an elementary $p$-group of rank $k$. 
In particular, we will focus on  the case where $G$ is abelian. In contrast to the other two types of ``large" group actions, here strategies balance a mix of knowledge of the quotient space with more general connectedness principles, which lead to structure theorems for topological invariants of the manifold. 
  \subsection{Positive Curvature and Large Symmetry Rank}

 Three fundamental results in this direction are the Maximal, Almost Maximal and Half-Maximal Symmetry Rank theorems due to Grove and Searle \cite{GrS1}, 
 Rong \cite{Ro} 
 and Fang and Rong \cite{FRo}, 
 and Wilking \cite{Wi2},
  respectively.  We present them together in one theorem.
 
 \begin{theorem}\label{sr} Let $M^n$ be a closed, positively curved manifold admitting an isometric and effective $T^k$-action. Then the following are true:
 
 \begin{enumerate}
 \item[1] \label{1} {\bf Maximal Symmetry Rank} (Grove and Searle
 \cite{GrS1})
If $k \geq \tfrac{n}{2}$, then $k=\lfloor \tfrac{n+1}{2}\rfloor$ and $M^n$ is  diffeomorphic to $S^{2k}$, $\mathbb{R}\mathrm{P}^{2k}$, $\CP^{k}$, or $S^{2k-1}/\zzz_l$ for some $l \geq 3$;
 \item[2] {\bf Almost Maximal Symmetry Rank}
   (Rong \cite{Ro}, and Fang and Rong \cite{FRo})
 \label{2} If $n\neq 6, 7$, $\pi_1(M)=0$, and $k= \lfloor\tfrac{n-1}{2}\rfloor$, then $M^n$ is homeomorphic to $S^{n}$, $\CP^{n/2}$, or $\HP^{2}$ for $k=3$ only; and
 \item[3] {\bf Half-Maximal Symmetry Rank}
  (Wilking \cite{Wi2})
 \label{3} If $n \geq10$,  $\pi_1(M)=0$, and $k \geq {\frac{n}{4}} + 1$, then $M^n$ is homeomorphic to $S^n$ or $\HP^{k-1}$ with $k = \frac{n}{4} + 1$, or homotopy equivalent to $\CP^{n/2}$.
 \end{enumerate}
 \end{theorem}
Observe first that the Maximal Symmetry Rank result can be improved to equivariant diffeomorphism with a linear $T^k$ action by work of Galaz-Garc\'ia \cite{G-G1}.
Additionally, by recent work of Kennard, Khalili Samani, and Searle \cite{KKSS}, The Half-Maximal Symmetry Rank result can be improved as follows: dropping the hypothesis of simple-connectivity, one can show that the only additional manifolds that occur are homotopy equivalent to $\RP^n$ and lens spaces of dimension $2k-1$, where $k=\frac{n}{4} + 1$.

Some comments on the proofs are in order, as many results stemming from these have leveraged the same techniques. 
The proof of  the Maximal Symmetry Rank result hinges on the fact that for the maximal symmetry rank, one can always find a circle subgroup of the $T^k$ acting fixed point homogeneously.  The Almost  Maximal Symmetry Rank result follows from  the Half-Maximal Symmetry Rank result and relies on Sullivan's homeomorphism classification of homotopy complex projective spaces and an analysis of the singular set of the group action  to improve the classification of  The Half-Maximal Symmetry Rank result to homeomorphism, as well as extend the result to dimensions $8$ and $9$. 
The proof of Half-Maximal Symmetry Rank result utilizes the theory of Error Correcting Codes, which give information about the dimensions of fixed point sets of involutions, and the following Connectedness Lemma of Wilking \cite{Wi2}.

\begin{theorem}[Wilking]\cite{Wi2} 
Let $M^n$ be a closed Riemannian manifold with positive sectional curvature. 
        If $N^{n-k}$ is a closed, totally geodesic submanifold of $M$, then the inclusion map $N^{n-k}\hookrightarrow M^n$ 
       is $(n-2k+1)$-connected.
\end{theorem}
 Recall that  a map $f: N\rightarrow M$ is called {\em $m$-connected}, if
$\pi_i(f):\pi_i(N)\rightarrow \pi_i(M)$ is an isomorphism for all $1\leq i\leq m-1$ and $\pi_m(f)$ is surjective. Using the Hurewicz isomorphism one can then make a similar statement about homology (and cohomology) groups. In Conclusion 3 of Theorem \ref{sr}, the Connectedness Lemma is leveraged by
combining Poincar\'e duality of $N$ with that of $M$, one obtains a certain periodicity  of the cohomology ring of  $M$, which in turn can be used to obtain homotopy equivalence.
These results have lead to significant progress on both Hopf Conjectures.
 As mentioned in the Introduction, for the first Hopf Conjecture, the seminal result of Hsiang and Kleiner established that in the presence of circle symmetry the conjecture is true.
 For the higher dimensional Hopf conjecture, which states that the product of two positively curved manifolds does not admit positive curvature, Amann and Kennard \cite{AKe} have shown that given $M$, a closed, simply-connected $n$-manifold, then the $2n$-dimensional product manifold $M\times M$ does not admit a metric of positive
 sectional curvature and an isometric torus action of rank $r> \log_{\frac{4}{3}}(2n-3)$.

For the second Hopf Conjecture, 
the main strategy to prove this conjecture with symmetries has been to use the fact, due to Kobayashi \cite{Ko}, 
that $$\chi(M)=\chi(M^{S^1}),$$ for some $S^1$ subgroup of an isometric $G$-action.
 P\"uttmann and Searle   \cite{PuS}
  and independently Rong    \cite{Ro}
   showed that for a $T^k$-action on a $2m$-dimensional manifold $M^{2m}$ with $k\geq \lfloor\frac{m-1}{2}\rfloor$ has $\chi(M)>0$. This lower bound for $k$ was quickly improved to $\lfloor \frac{m-2}{4}\rfloor$ by Rong and Su \cite{RoSu}, 
   and to $\lfloor \frac{m}{5}\rfloor$ by Su and Wang \cite{SuWa}.
    Using the method of Steenrod squares, Kennard \cite{Ke} 
    was able to improve this lower bound to $\log_{2}(2m-2)$. More recently Kennard, Wiemeler, and Wilking \cite{KeWieWi} found a general lower bound for $k$, completely independent of dimension. \begin{theorem}[Kennard, Wiemeler, and Wilking \cite{KeWieWi}] If $M^{2m}$ is an even-dimensional, connected, closed, positively curved Riemannian manifold whose isometry group has rank at least five, then $\chi(M) > 0$.
 \end{theorem}
 It is reasonable to expect that the tools used in the proof of this theorem will be useful in many different contexts.
 Finally, Nienhaus \cite{Ni} has announced an improvement on this theorem, claiming to be able to lower the bound to $4$.
 
 \subsection{Non-negative Curvature and Large Symmetry Rank}

In strong contrast to the positive curvature case and in part due to the lack of guaranteed fixed point sets of torus actions, the maximal symmetry rank of a closed, simply-connected, non-negatively curved manifold has not yet been established in all dimensions.  
\begin{msrconj} Let $T^k$ act isometrically and effectively on
$M^n$, a closed, simply connected, non-negatively curved Riemannian manifold. Then  
$k\leq \lfloor 2n/3\rfloor$ and 
when $k= \lfloor 2n/3\rfloor$, $M^n$ is equivariantly diffeomorphic to  
$\mathcal{Z}/T^m$ with a linear $T^k$-action,  where 
$$\mathcal{Z}=  \prod_{i\leq r} S^{2n_i-1} \times\prod_{i>r} S^{2n_i},$$
 with  $n_i\geq 2, \,\, r= 2\lfloor 2n/3\rfloor-n, \,\,0 \leq m \leq 2n \mod 3,$
and the $T^m$-action on $\mathcal{Z}$ is  free and linear.
\end{msrconj} 
We can then ask the following question.
 \begin{problem}\label{msr} Let $\mathcal{M}_0^n$  denote the class of closed, simply-connected, non-negatively curved $n$-manifolds. For all $M\in \mathcal{M}_0^n$ establish the upper bound for the maximal symmetry rank  and classify all $M\in \mathcal{M}_0^n$ 
of maximal and almost maximal symmetry rank.
\end{problem}

Returning to the Maximal Symmetry Rank Conjecture, 
work of Galaz-Garc\'ia and Searle \cite{G-GS1}, 
Galaz-Garc\'ia and Kerin \cite{G-GK}, 
and of Escher and Searle \cite{ES1} shows that the conjecture holds in dimensions $n\leq 9$. The upper bound for dimensions less than or equal to $12$ was established by Galaz-Garc\'ia and Searle    \cite{G-GS1} 
and by Escher and Searle \cite{ES1}.

\begin{figure}[!htb]
\centering
\vspace{5cm}
\hspace{-8.5cm}
\begin{picture}(0,0)
\includegraphics[scale=0.48]{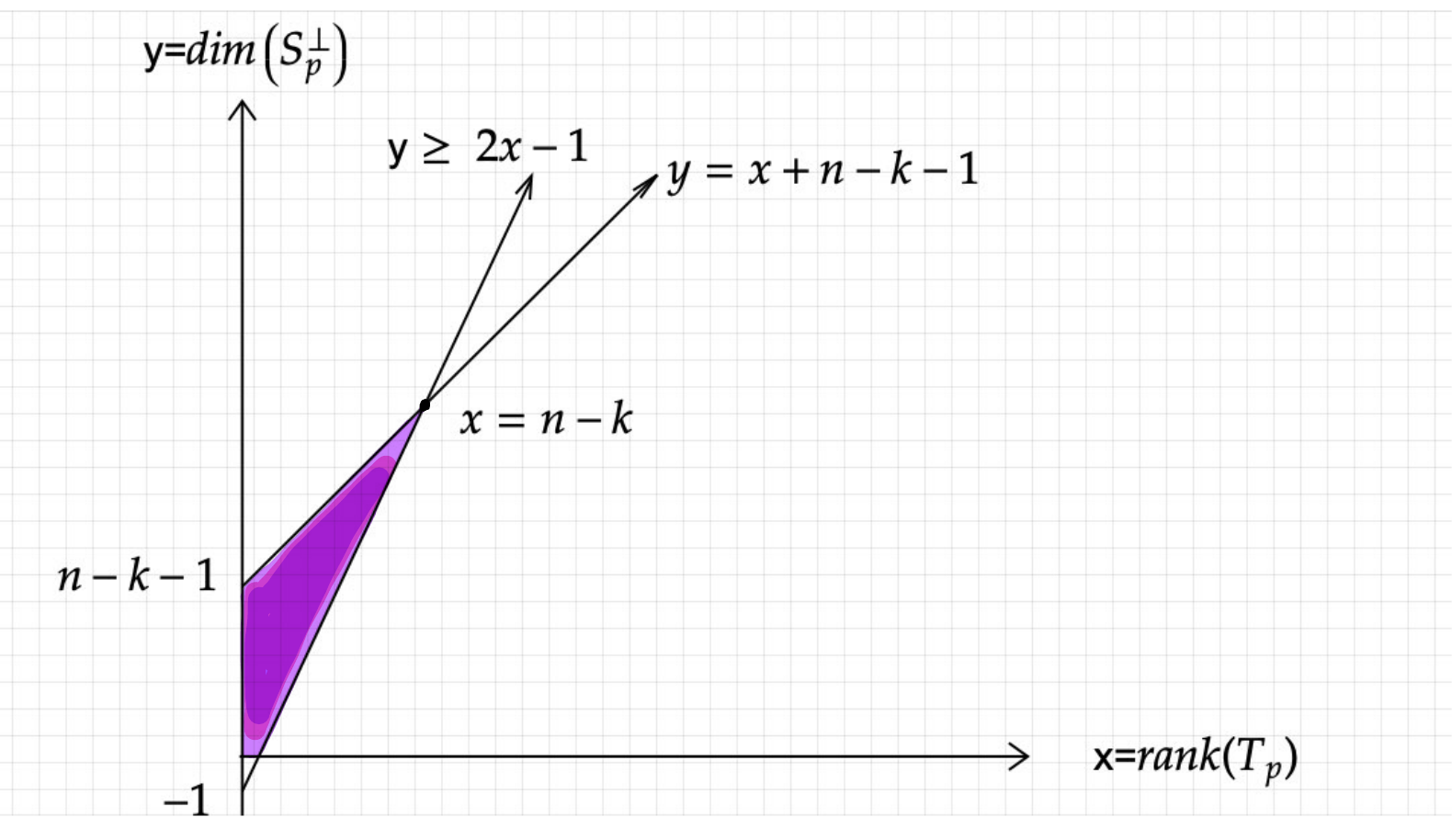}
\end{picture}%
\setlength{\unitlength}{3947sp}%
\begingroup\makeatletter\ifx\SetFigFont\undefined%
\gdef\SetFigFont#1#2#3#4#5{%
  \reset@font\fontsize{#1}{#2pt}%
  \fontfamily{#3}\fontseries{#4}\fontshape{#5}%
  \selectfont}%
\fi\endgroup%
\caption{The Range of Possible Ranks for the Isotropy Subgroup of a Torus Action}
\label{msrimage}
\end{figure}

We observe that $n-k$ is the maximal possible rank of an isotropy group of a torus action: the dimension of the unit normal sphere to $p\in T(p)$ is a function of the rank of $T_p$ and 
$\rk(T_p)\leq \lfloor (\dim(S^{\perp}_p)+1)/2\rfloor$ by the Maximal Symmetry Rank Theorem, see Figure \ref{msrimage}.
Notably, the $n$-manifolds described in Part 2 of the Maximal Symmetry Rank Conjecture all admit  $T^k$-actions such that  there is a point 
$x \in M^n$ for which $\rk(T_x)=n-k$, or equivalently, $\dim(T(x))=2k-n$. We call such actions  {\em isotropy-maximal}.

In particular, for an isotropy-maximal $T^k$-action, $M^{T^{n-k}}\neq \emptyset$. In fact,  any component of $M^{T^{n-k}}$ is contained in a {\em generalized characteristic submanifold} of $M$. That is, there is some circle subgroup $S^1\subset T^k$ with a codimension two fixed point set component, $F$, that contains a $T^{n-k}$-fixed point set component. When $M^{T^{n-k}}$ is $0$-dimensional, $F$ is simply called a characteristic submanifold.
In particular, one sees that an isotropy-maximal $T^k$-action is an example of a nested $S^1$-fixed point homogeneous action, generating a tower of 
nested fixed point sets of subtori of the $T^k$-action in $M$.

The case where $\dim(M^{T^{n-k}})=0$ corresponds to that of {\em torus manifolds}.   Wiemeler   in    \cite{Wie}
 classified closed, simply-connected, non-negatively curved torus manifolds, finding that they are all   equivariantly diffeomorphic to a quotient of a free linear torus action of $\mathcal{Z}$, as in the Maximal Symmetry Rank Conjecture.

Escher and Searle in \cite{ES1} generalized this result to all isotropy-maximal torus actions on closed, simply-connected, non-negatively curved $n$-manifolds, showing that they are all  equivariantly diffeomorphic to a quotient of a free linear torus action of 
$\mathcal{Z}$, as in the Maximal Symmetry Rank Conjecture.
Indeed,  if the Bott Conjecture holds, one can combine the isotropy-maximal classification result of Escher and Searle \cite{ES1} with the work of Galaz-Garc\'ia, Kerin, and Radeschi \cite{G-GKR} mentioned in Remark \ref{bottGGKR}, to show that the Maximal Symmetry Rank Conjecture holds. 

Since then, Dong, Escher, and Searle in \cite{DES} have extended the result to almost isotropy-maximal torus actions, where an action is {\em almost isotropy-maximal} if 
there is a point $x \in M$ such that 
the dimension of its isotropy group is $n-k-1$,  or, equivalently
there is a point $x \in M$ whose orbit is of dimension   $2k-n+1$. In particular, the manifolds obtained are as in the isotropy-maximal classification.

We observe that $n$-manifolds of almost maximal symmetry rank have also been classified in dimensions less than or equal to $6$ by independent work of Kleiner   \cite{Kl} 
and of Searle and Yang   \cite{SY} 
in dimension $4$, by work of Galaz-Garc\'ia and Searle   \cite{G-GS2} 
in dimension $5$ and by work of Escher and Searle \cite{ES2} in dimension $6$. 
Notably, it is only in dimension $5$ that we observe any difference with the maximal symmetry rank classification, as the Wu manifold, $SU(3)/SO(3)$, which is not the quotient of a linear torus action on a product of spheres of dimensions greater than or equal to three, appears. 
The work above suggests that one approach to Problem \ref{msr} would be to begin by classifying torus actions on closed, simply-connected, non-negatively curved manifolds via the rank of the largest possible isotropy group, beginning with those that have rank $n-k-2$. 

Remark that by work of B\"ohm and Wilking \cite{BoWi}, 
a closed, simply-connected, non-negatively curved manifold admits a metric of positive Ricci curvature.
Thus, another approach  when studying the class $\mathcal{M}_0^n$ is to consider the {\em $k$th-intermediate Ricci curvature}, a curvature interpolating between  sectional curvature and Ricci curvature.  Of natural interest  is the case where this curvature is positive. 
 We say an $n$-dimensional Riemannian manifold $(M^n,g)$ has {\em positive $k$th-intermediate Ricci curvature} for $k\in\{1, \hdots,n-1\}$ if, for any choice of orthonormal vectors $\{u, e_1, \hdots, e_k\}$, the sum of sectional curvatures $\Sigma_{i=1}^k \sec(u, e_i)$ is positive. 
Observe that $\Ric_1 > 0$ is equivalent to positive sectional curvature, and $\Ric_{n-1} > 0$ is equivalent to positive Ricci curvature. Furthermore, if $\Ric_k > 0$, then $\Ric_l > 0$ for all $l \geq k$.

One advantage to studying such curvature lower bounds is that many of the results for positive sectional curvature, such as the existence of fixed points and Wilking's Connectedness Lemma extend to this curvature bound. 
A natural starting point to consider symmetries in the presence of positive intermediate Ricci curvature is to consider those with 
$\Ric_2>0$. Mouill\'e \cite{Mou} in recent work has been able to extend the Maximal Symmetry Rank Theorem  to closed, $n$-manifolds with $\Ric_2>0$, showing that the symmetry rank of such manifolds is the same as for the positive curvature case, and moreover obtains a similar classification result.

Finally, while a product of spheres $\prod_i^m S^{n_i}$ with the product metric has $\Ric_k>0$ for $k\geq \max_{i\in\{1, \hdots, m\}}\{ 1+ \sum_{j\neq i} n_j\}$, there are also metrics of positive $\Ric_{k'}$ with $k'<k$ on such products. In particular, in Example 2.3 of Mouill\'e \cite{Mou}, he shows that one can put $\Ric_2>0$ metrics on $M^6=S^3\times S^3$ with $T^3$ symmetry, as well as $\Ric_2>0$ on quotients by free torus  actions $M^6$. These results are consistent with the almost maximal symmetry rank classification results mentioned earlier and suggest a connection between closed, simply-connected manifolds of non-negative curvature with large symmetry rank and  $Ric_k>0$ for small $k$. 
It would be interesting to explore this relationship, if any exists.

\subsection{Positive and Non-negative Curvature and Large Discrete Symmetry Rank}
  
 As mentioned earlier, a $\zzz_p^k$-action need not 
 have fixed points. To guarantee the existence of fixed points, one may, for example, assume that $p$ is larger than $\mathcal{C}_n$, Gromov's estimate for the total Betti number of a closed $n$-manifold.  
  Work of Yang   \cite{Y} 
  and Hicks   \cite{Hi} 
  showed that for closed, simply-connected, non-negatively curved $4$-manifolds with discrete $p$-symmetry rank $1$ and $2$, with $p$ suitably chosen, one can bound the total Betti number of $M^4$ above by $7$ and $5$, respectively.   In higher dimensions,  Fang and Rong   \cite{FRo} 
  extend the Maximal Symmetry Rank Theorem as follows.
  \begin{theorem}[Fang and Rong]  \cite{FaRo1} 
  Let $M^n$ be a closed, simply-connected, positively curved $n$-manifold admitting an isometric and effective $G$-action, with $G=\zzz_p^k$,  $p>\mathcal{C}_n$, and $k\geq \lfloor \frac{n}{2} \rfloor$. Then the following hold:
  \begin{enumerate}
  \item If $n=2m$,   then $k=m$ and $M^n$ is homeomorphic to $S^n$ or $\CP^m$.
  \item If $n=2m+1$ and the $G$ action is extended to an  isometric  and effective  $S^1\times G$-action, then $k=m$ and $M^n$ is homeomorphic to $S^n$.
  \end{enumerate}
  \end{theorem}

Fang and Rong   \cite{FaRo1} 
also obtain an analog of the Half Maximal Symmetry Rank Theorem for  discrete $p$-symmetry rank bounded below by approximately $3n/4$. Further reductions of the discrete $p$-symmetry rank for $p>\mathcal{C}_n$ have resulted in confirmation of the second Hopf Conjecture for this class of manifolds by Rong and Su \cite{RoSu}, 
as well as work by Wang \cite{Wa}, 
extending work of Rong   \cite{Ro2} 
showing that  $\pi_1(M^n)$ is cyclic provided $\symrk(M^n)>\frac{n}{4}+1$, to $T^1\times \zzz_p^k$-actions with  $k>\frac{n+1}{4}$.

More recently, Kennard, Khalili Samani, and Searle in \cite{KKSS} consider positively curved Riemannian manifolds admitting an isometric action of $\zzz_2^r$   with a fixed point. 
They extend the Maximal Symmetry Rank results of   Grove and Searle    \cite{GrS1} 
and the Half-Maximal Symmetry Rank results of Wilking   \cite{Wi2} 
to obtain the following three results, noting that the first result is an easy consequence of the work of Fang and Grove in \cite{FGr}.

\begin{theorem}[Kennard, Khalili Samani, and Searle \cite{KKSS}]\label{thm:n}
Let $M^n$ be a closed, positively curved manifold such that $\zzz_2^r$ acts isometrically on $M$ with $x\in M^{\zzz_2^r}$.
Then the following hold:
\begin{enumerate}
\item\label{11}  If $r \geq n$, then $r=n$ and  $M$ is equivariantly diffeomorphic to $S^r$ or $\R{\mathrm P}^r$  with a linear $\zzz_2^r$-action.
\item\label{21} If $n\geq 24$ and
$r \geq \tfrac{n+1}{2}$, 
then 
$M$ is homeomorphic to $S^n$, or homotopy equivalent to  $\rrr{\mathrm P}^n$, 
$\ccc{\mathrm P}^{r-1}$, or $S^{2r-1}/\zzz_k$ for some $k \geq 3$ and $r=\lceil \tfrac{n+1}{2}\rceil$.
\item\label{31} If $n \geq 15$ and  $r \geq \tfrac{n+3}{4} \,+ \,1$. 
Then at least one of the following occurs:
	\begin{enumerate}
	\item For any subgroup of $\zzz_2^r$ with corank at most four, the fixed point set component $F^m$ at $x$ is homotopy equivalent to $S^m$, $\rrr{\mathrm P}^m$, $\ccc{\mathrm P}^{\frac m 2}$, or $S^m/\zzz_k$ for some $k \geq 3$; or
	\item $M^n$ is a simply connected integer cohomology $\hh{\mathrm P}^{r-2}$ and $r=\tfrac{n}{4} +2$.
	\end{enumerate}
	\end{enumerate}
\end{theorem}

In Figure \ref{bothmsr}  below, the  Torus and $\zzz_2$-Torus Symmetry Rank results are displayed graphically for closed, simply-connected, positively curved manifolds.
\begin{figure}[!htb]
 \centering
 {%
      \includegraphics[width=0.474\textwidth]{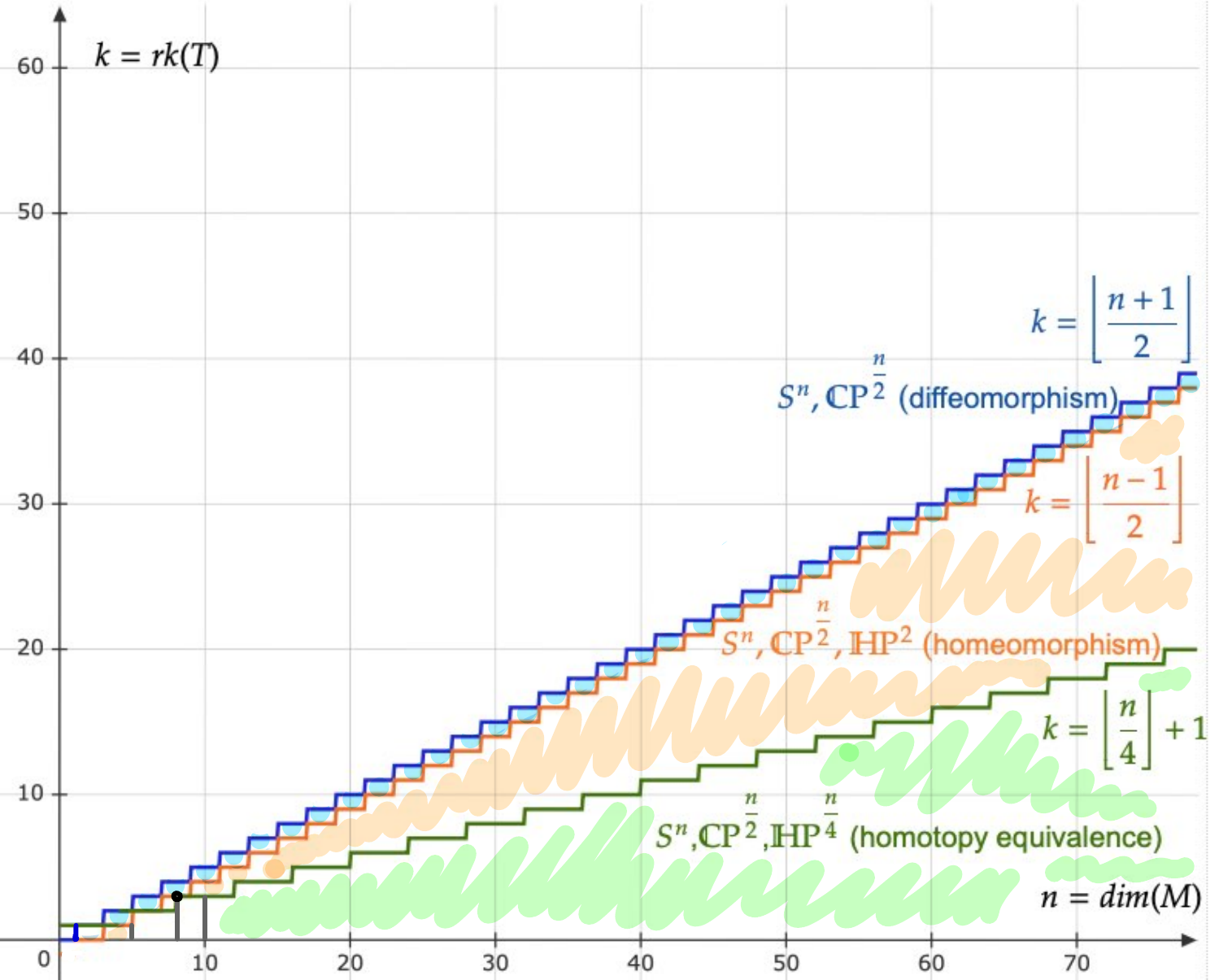}}
      \label{fig:image-a}
 \qquad
 {%
      \includegraphics[width=0.448\textwidth]{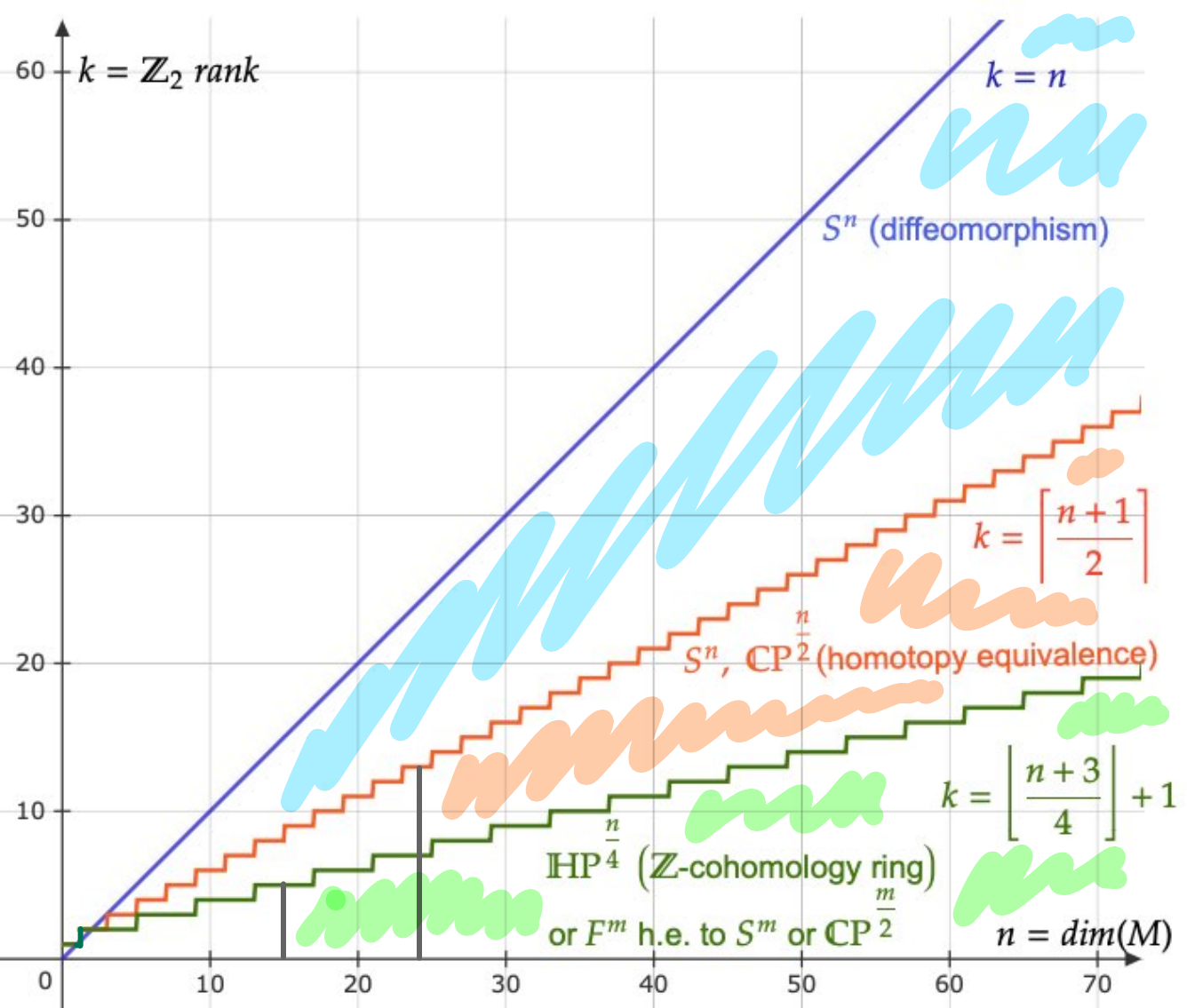}}
      \label{fig:image-b}
      \caption{The Torus and $\zzz_2$-Torus Symmetry Rank Theorems for simply-connected $M$}
       \label{bothmsr}
\label{msrt}
\end{figure}
The proof of Conclusion \ref{11} relies on understanding the quotient space of $M$ by an involution fixing a codimension $1$ submanifold. On the other hand, the techniques used to prove Conclusion \ref{21} include Wilking's Connectedness Lemma, and refined estimates for Error Correcting Codes, and an inductive argument. Conclusion \ref{31} however, is not proven by induction and relies heavily on the Borel Formula.  

The following question, albeit ambitious, is natural, given the path taken with torus actions to achieve the recent result of Kennard, Wiemeler, and Wilking.
\begin{question}
How much of the work in Kennard, Wiemeler, and Wilking \cite{KeWieWi} can be adapted to the case of $\zzz_2^k$-actions with $M^{\zzz_2^k}\neq\emptyset$. \end{question}

Since $\chi(M)\equiv \chi(M^{\zzz_2^k})\mod 2$, one does not expect to achieve an analog of their Euler characteristic result, however, other important results in Kennard, Wiemeler, and Wilking \cite{KeWieWi} may yield analogs for the $\zzz_2$ case.

 \section{Acknowledgements} The author is extremely grateful to Christine Escher and Fred Wilhelm for their willingness to read multiple draft versions of this article and for their valuable feedback. She also thanks the referees for their careful reading and suggestions for improvement. The author's work is partially supported by NSF-DMS grant \#2204324.


\end{document}